\documentclass{article}
\usepackage{a4wide}
\usepackage{amssymb}
\usepackage{amsfonts}
\usepackage{amsxtra}
\usepackage{latexsym}

\begin{document}

\newtheorem{theo}{Theorem}[section]
\newtheorem{theo*}{Theorem}
\newtheorem{lemma}[theo]{Lemma}
\newtheorem{definition}[theo]{Definition}
\newtheorem{remark}[theo]{Remark}
\newtheorem{cor}[theo]{Corollary}
\newtheorem{prop}[theo]{Proposition}
\newtheorem{example}[theo]{Example}

\newcommand{\N}{{\mathbb{N}}}
\newcommand{\Z}{{\mathbb{Z}}}   
\newcommand{\Q}{{\mathbb{Q}}}
\newcommand{\R}{{\mathbb{R}}}
\newcommand{\C}{{\mathbb{C}}}
\newcommand{\M}{{\mathbb{M}}}
\newcommand{\F}{{\mathbb{F}}}
\newcommand{\LL}{{\mathbb{L}}}

\newcommand{\g}[1]{{\mathfrak {#1}}}
\newcommand{\cc}[1]{{\cal {#1}}}
\newcommand{\ul}{\underline }
\newcommand{\ol}{\overline }

\newcommand{\qed}{\hspace*{\fill}$\Box$}
\newcommand{\Qed}{\hspace*{\fill}$\Box \Box$}

\newcommand{\md}{{\rm mod}\,\,}
\newcommand{\ee}{{\rm zip}\,}
\newcommand{\ei}{{\rm zap}\,}


\title{Self-similar carpets over finite fields}
\author{Mihai Prunescu \thanks{ Brain Products, Freiburg, Germany, and Institute of Mathematics of the 
Romanian Academy, Bucharest, Romania.
mihai.prunescu@math.uni-freiburg.de.}}
\date{}
\maketitle

\begin{abstract}

\parindent 0 cm 

In \cite{PL} an informal algorithm 'to display interesting numeric patterns' is described 
without any proof. We generalize this algorithm over arbitrary finite fields of characteristic $p$
and we prove that it really generates self-similar carperts, provided that
they contain at least one zero in the first $(p+1)/2$ lines. 
For the fields $\F_p$ the generalized 
algorithm produces $p-1$ different self-similar carpets. 
These self-similar carpets are classified according to
their arithmetic and their groups of symmetry. All this phenomena can be also interpreted in the framework
of the aperiodic tilings.

\thanks{A.M.S.-Classification: 11A07, 28A80.}
\end{abstract} 


\section{Introduction}

In \cite{PL} is presented an informal algorithm 'to produce interesting numerical patterns'. 
Let $n$ be a fixed natural numbers $ > 2$. One takes a rectangular matrix,
completes the first row and the first column with ones, and recursively computes the other elements as 
$( N + NW + W ) \mod n$,
where $N$, $NW$ and $W$ are the neighbors in the corresponding directions. Finally, one can produce an 
image following a fixed correspondence of the rests modulo $n$ and a list of colours. The autors observe and state 
that for primes $n = p$ the patterns are self-similar, but don't prove this. For the notion of 
self-similarity they citate Mandelbrot's monograph
\cite{BM}. In \cite{BM}, Chapter 14,  there is a hint to a similar construction of the 
original Sierpinski Carpet atributed to Rose (see \cite{Ro}). One referee kindly informed us  about the paper \cite{PL1}
where the authors introduced a generalization of this generation rule $( N + m \cdot NW + W ) \mod n$ for a constant 
$m \in \N$ and commented on the associated generalized Fibonacci sequences, but didn't do any graphical interpretation.

One goal of this paper is to prove the conjecture suggested in \cite{PL}: if $n = p$ odd prime, the rescaled images 
produced by the a recurent rule like in \cite{PL1} converge to a symmetric self-similar set. For $m=0$ one recovers 
the known case of  Pascal's Triangle, with a group of symmetries consisting of two elements. 
For $m=-1$ one gets the full square - which is self-similar but boaring. For
$m=1$ one gets self-similar sets of \cite{PL} that have as group of symmetries the full dihedral group $D_8$. 
For the other values of $m$ one gets new self-similar patterns with group of symmetries isomorphic with Klein's group $K_4$. 
This paper is a generalization of author's paper \cite{Pr} where only the case $m=1$ was treated. 

We also prove that those carpets over arbitrary finite fields, which have zeros, are self-similar and respect a
simple law of symmetry. 

The results can be also understood as existence of aperiodic tilings of a quadrant of the plane

\begin{definition}\rm Let $\F_q$ be an arbitrary finite field and fix an element $m \in \F_q$. 
The matrices occurring in this article 
are always indexed from $0$ and have elements in 
$\F_q$, if not otherwise specified. Let the prime $p$ be the characteristic of the finite field, $q = p^k$ for some $k$.
Let $M_d = (a_{i,j})$ be the $p^d \times p^d$ 
matrix constructed following the recurence
$a_{i,0} = a_{0,j} = 1$ and $ a_{i,j} = a_{i-1,j} + m \cdot a_{i-1,j-1} + a_{i,j-1}$. 
The matrix $M_1$ shall be denoted by $F(p,m)$ and shall be called
{\bf fundamental block}. 
\end{definition}

\begin{definition}\rm The black and white image $I^d$ is defined as follows: 
one tiles the compact square $[0,1] \times [0,1]$ in $p^d \times p^d$
many equal squares $S_{i,j}$, and excludes the interior of $S_{i,j}$ if and only if 
$a_{i,j}=0$. 
\end{definition}

\begin{definition}\rm
The self-similar set in question shall be 
$I = \lim I^d$. The name of the variable $d$ is chosen to mean the {\bf depth} of 
the recursive approximation of $I$. The limit operator
can be understood in the sense of the Hausdorff metric for compact subsets of $\R^2$. 
\end{definition}

\begin{definition}\rm The coloured image of $M_d$ is defined assuming a fixed correspondence
between the rests modulo $p$ and a family of colours. In order to make this consistent with the
black and white image, suppose that the colour {\it white} always corresponds to zero.
\end{definition}


\section{The recurrent function}

In this section $K$ is the notation of an arbitrary field.

{\bf Definition}: Let $m \in K$ be fixed. 
We consider the function $f : \N \times \N \rightarrow K$ recursively defined by the conditions 
$f(n,0)=f(0,k)=1$ and:
\[f(n,k) = f(n,k-1)+ m \cdot f(n-1,k-1)+f(n-1,k)\]
for $n,k\geq 1$.

\begin{lemma}\label{rec}
The function $f$ is symmetric and satisfies:
\[f(n,k) = \sum \limits _{a=0} ^{\min (n,k)} m^a \binom{n}{a} \binom{n+k-a}{k-a}.\]
\end{lemma}

{\bf Proof}: The symmetry follows from the symmetry of the recurrence formula and of the initial conditions.
To compute $f$, use the method of generating functions, see \cite{Wf}. Define the generating function 
$A_n(x) = \sum \limits _{k\geq 0}f(n,k) x^k$. It follows:
\[A_{n+1}(x) = \sum \limits _{k\geq 0} f(n+1,k) x^k = 1 + \sum \limits _{k \geq 1} \Big{(}f(n,k) + f(n+1,k-1) 
+ m \cdot f(n,k-1)\Big{)}x^k =\]
\[=\Big{(}1+\sum \limits _{k\geq 1} f(n,k)x^k\Big{)} +x\sum \limits _{k\geq 0} 
f(n+1,k) x^k  + m x \sum \limits _{k\geq 0} f(n,k) x^k =
A_n(x) + xA_{n+1}(x) + m x A_n(x).\]
This recurrence have the solution:
\[A_n(x) = {\Big{(} \frac{1}{1-x}\Big{)} }^{n+1} {(1+mx)}^n.\]
Using that ${(1+mx)}^n = \sum \limits _{k \geq 0} \binom{n}{k}m^kx^k$ and that $ { \Big{(} \frac{1}{1-x} \Big{)} }^{n+1} 
= \sum \limits _{k \geq 0} \binom{n+k}{k}x^k$, one gets the Lemma. \qed

One remarks that the terms
\[t(a,k,n) = m^a \binom{n}{a} \binom{n+k-a}{k-a} = m^a \frac{(n+k-a)!}{a!(k-a)!(n-a)!}.\]
are itself symmetric in $n$ and $k$.


\section{Tensor powers and the automorphism of Frobenius}

In this section we prove some properties of the fundamental block $F(p,m) \in {\cal M}_{p \times p}(\F_p)$. 
Recall the notation $F(p,m) = (a_{i,j})$
with $i$ and $j = 0, \dots , p-1$. From Lemma \ref{rec} it follows already that $F(p,m)$ is a symmetric matrix, 
like all other $M_d$.

\begin{lemma}\label{lastrow} 
The last column and the last row of $F(p,m)$ are exactly:
\[1, \, -m, \, (-m)^2, \, \dots ,\, (-m)^{p-1}.\]
This works also for $m=0$.
\end{lemma}

{\bf Proof}: Take $k \leq n = p-1$ and work over $\F_q$. For $a < k$ the term:
\[t(a,p-1,k) = m^a \binom{p-1}{a} \binom{p-1+k-a}{k-a} =
m^a \binom{p-1}{a}\cdot p \cdot \dots =0,\]
so all these terms do not contribute in $\F_q$. For the last term one has:
\[t(k,k,p-1) =  m^k\frac{(p-1)\dots(p-k)}{k!} = m^k {(-1)}^k \frac{k!}{k!} = (-m)^k. \]

\qed

\begin{definition}\rm Let $R$ be some commutative ring and $A = (a_{i,j})\in {\cal M}_{s \times t}(R)$, 
$B \in {\cal M}_{u \times v}(R)$ 
two matrices. Then the {\bf tensor product} $A \otimes B$ is a matrix in ${\cal M}_{su \times tv}(R)$ having the 
block-representation $(a_{i,j}B)$.  If $A_1$, $A_2$, $\dots$, $A_n$ are arbitrary matrices, we denote the tensor term:
\[(( \dots ((A_1 \otimes A_2) \otimes A_3) \dots ) \otimes A_{n-1}) \otimes A_n .\]
by:
\[A_1 \otimes A_2 \otimes \dots \otimes A_{n-1} \otimes A_n.\]
For all $n \geq 1$ we define the {\bf tensor power} $A^{\otimes n}$ of $A$ 
inductively by: $A^{\otimes 1} = A$ and $A^{\otimes (n+1)} = A^{\otimes n} \otimes A$.
\end{definition}

The following remark expresses the {\it principle of substitution} used in constructing self-similar sets.

\begin{remark}\label{tenssimilar}
For some $n \geq 2$ consider a matrix $A \in {\cal M} _{n \times n}(\{0,1\})$ 
containing at least one zero and at least two ones.
Let $I^d$ be the black and white image associated to  $A^{\otimes d}$. 
Then $I^d$ is the $d$-th step in the transfinite construction of a non-trivial self-similar
set $I = \lim I^d$.
\end{remark}

Now is the time that the automorphism of Frobenius enters the scene. 

\begin{definition}\rm
The automorphismus of Frobenius $\varphi : \F_q \rightarrow \F_q$ is defined by 
$\varphi(x) = x^p$. This automorphism generates the Galois group $G(\F_q / \F_p)$.
\end{definition}

\begin{lemma}\label{cells}
Let $F = F(p,m)$ be a fundamental block for some $m \in \F_q$. Consider the matrix in construction:
\[
\begin{array}{cc}
\alpha F & \beta F \\
\gamma F & \cdot
\end{array}
\]
with $\alpha$, $\beta$, $\gamma$ $\in \F_q$. By application of the recurrent rule  one gets: 
\[
\begin{array}{cc}
\alpha F & \beta F \\
\gamma F & \delta F
\end{array}
\]
with $\delta = \varphi(m)\alpha +\beta + \gamma$.
\end{lemma}

{\bf Proof}: Denote $-m$ with $\lambda$.
There is only one element $x$ where one can start to apply the recurrent rule: 

\[
\begin{array}{cccccc}
\dots&\dots& \lambda ^{p-3}\alpha & \beta&\dots&\dots\\
\dots&\dots&\lambda ^{p-2}\alpha& \beta&\dots&\dots\\
 \lambda ^{p-3}\alpha&\lambda ^{p-2}\alpha& \lambda ^{p-1}\alpha& \beta&\lambda \beta& \lambda ^2 \beta\\
 \gamma& \gamma& \gamma & x&\cdot&\cdot \\
\dots&\dots&\lambda \gamma&\cdot&\cdot&\cdot \\
\dots&\dots&\lambda ^2 \gamma&\cdot&\cdot&\cdot
\end{array}
\]
Applying the recurent rule along the first row and along the first column to be completed yelds:
\[
\begin{array}{cccccc}
\dots&\dots& \lambda ^{p-3}\alpha& \beta&\dots&\dots\\
\dots&\dots&\lambda ^{p-2}\alpha& \beta&\dots&\dots\\
 \lambda ^{p-3}a&\lambda ^{p-2}\alpha& \lambda ^{p-1}\alpha& \beta&\lambda \beta& \lambda ^2 \beta\\
 \gamma& \gamma& \gamma& \delta& \delta&\delta \\
\dots&\dots&\lambda \gamma& \delta&\cdot&\cdot \\
\dots&\dots&\lambda ^2 \gamma& \delta &\cdot&\cdot
\end{array}
\]
where $\delta =(-m)^{p-1}\alpha * m+\beta+\gamma = m^p \alpha + \beta + \gamma =$ $
 \varphi(m) \alpha + \beta + \gamma$ in $\F_q$. 
The recurrent rule is linear, so a constant $\delta$ row together with a constant $\delta$ 
column generate $\delta F$. \qed

\begin{definition}\rm
For a matrix $A=( a_{i,j})$ over $\F_q$, let $\varphi (A)$ be the matrix $(\varphi (a_{i,j}))$.
\end{definition}

\begin{theo}\label{tensorpower} Recall that $M_d$ is the $p^d \times p^d$ 
matrix computed by the recurrent rule over the finite field $\F_q$ and $F = F(p,m) = M_1$ is 
the fundamental block. Then for all $d \geq 1$:
\[M_d = \varphi ^{d-1} (F) \otimes \varphi ^{d-2}(F) \otimes\dots \otimes \varphi(F) \otimes F. \]
\end{theo}

{\bf Proof}: The proof works by induction. For $d=1$ this is true by definition. 
Suppose that $M_d$ fulfills the stated relation and consider $M_{d+1}$. 
Being computed by the same recurent rule, the $p\times p$ left upper block of $M_{d+1}$ is a copy of $F$. 
Applying Lemma \ref{cells} for $(a,b,c) = (0,0,1)$ or $(0,1,0)$ one gets that a copy of $F$ continued by a 
row of ones $F^{111\dots 1}$ horizontally generates copies of $F$ like $FFFF \dots F$ and that this happens also 
vertically if the first column of $F$ is downwards extended with ones. Thus in the block-representation of $M_{d+1}$ 
with $p\times p$ blocks, the first line and the first column consist of copied fundamental blocks:
\[
\begin{array}{ccccc}
 F& F& F& F& \dots  \\
 F&b_{1,1}F&b_{1,2}F&b_{1,3}F& \dots  \\
 F&b_{2,1}F& \cdot & \cdot & \dots  \\
 F&b_{3,1}F& \cdot &\cdot  & \dots  \\
 \vdots &\vdots  &\vdots&\vdots&  \ddots 
\end{array}
\]
Here the $x_{i,j}$ are such that all $x_{k,0} = 1$, $x_{0,n} = 1$ and $x_{i+1,j+1} = \varphi(m) x_{i,j} + x_{i+1,j}
+x_{i,j+1}$. 
One gets that $M_{d+1} = X \otimes F$ where the matrix $X$
over $\F_q$: (i) has $p^d \times p^d$ elements, (ii)
has elements $x_{i,j}$ as above. This means that 
$X = M_d$ for the process made with the coefficient $m' = \varphi(m)$, so by induction, 
denoting $F(p, \varphi(m)) = F'$:
\[X = \varphi ^{d-1} (F') \otimes \varphi ^{d-2}(F') \otimes\dots \otimes \varphi(F') \otimes F'. \]
If we substitute $F' = F(p, \varphi(m)) = \varphi (F(p,m)) = \varphi (F)$, one gets:
\[X = \varphi ^d (F) \otimes \varphi ^{d-1}(F) \otimes\dots \otimes \varphi ^2(F) \otimes \varphi (F). \]
Substituting this $X$ in $M_{d+1} = X \otimes F$ one gets the desired relation. \qed

\begin{cor}\label{fractalInGeneral}
For all finite fields $\F_q$ and all $m \in \F_q$, if the fundamental block $F(p,m)$ contains 
at least a zero,  the black and white image $I^d$ of $M_d$ is the $d$-th 
step in the transfinite construction of a non-trivial self-similar set $I$. 
\end{cor}

{\bf Proof}: For any matrix $A$ over $\F_q$, let $\delta( A)$ be the matrix obtained 
by substituting every non-zero element with one.
Let $\iota(B)$ the black and white image of the matrix $B$. Then:
\[ I^d = \iota \delta (M_d) = \iota \delta (\bigotimes \limits _{i=d-1} ^0 \varphi ^i(F)) = \iota (D^{\otimes d}),\]
where $D = \delta(F)$ is a $\{0,1\}$-matrix, $F=F(p,m)$ is the fundamental block, and $\varphi$ is Frobenius' automorphism 
extended for matrices. Now the principle of substitution works. \qed

\begin{lemma}\label{existencezero}
If $m \in \F_p$,
the fundamental block $F(p,m)$ contains zeros if and only if $m \neq -1$. In this situation it contains in the row 
$i=1$ exactly one zero: 
\[a_{1,k} = 0 \,\, \leftrightarrow \,\, \F_p \,\models \, k = -(m+1)^{-1}.\]
Note: in general there are many other zeros in the fundamental block.

\end{lemma}

{\bf Proof}: The element $a_{1,k} = km + (k+1) = k(m+1) +1 $ which is zero only for $k = -(m+1)^{-1}$, existing for all
$m \neq -1$ in $\F_p$. Every such $k$ has a representative between $1$ and $p-1$ inclusively.
If $m=-1$ the matrix $F(p,-1)$ contains only the repeated element $1$. 
\qed

Now from Remark \ref{tenssimilar} and from the Lemma \ref{tensorpower} the main result follows:

\begin{cor}\label{fractalOverFp}
For all primes $p$ and all $m \in \F_p \setminus \{-1\}$ the black and white image $I^d$ of $M_d$ is the $d$-th 
step in the transfinite construction of a non-trivial self-similar set $I$. For $m=-1$ the set $I$ is the full
square. 
\end{cor}

It is worth to point out that:

\begin{cor}
The Pascal Triangle modulo $p$ (got for $m=0$) and the Passoja-Lakhtakia Carpets (got for odd $p$ and $m=1$) 
are non-trivial self-similar sets.
\end{cor}

\begin{example}\rm

Differently as for the black and white images, the coloured images are only symmetric according to the first diagonal, 
excepting for the colour 'white' (the holes).
However, strange and optically exciting coloured patterns arise.
The following example shows the step $M_2$ for $p = 3$ and $m=1$, 
a step in the construction for the celebrated Sierpinski Carpet, used in \cite{S}. The zeros are not displayed.

\[
\arraycolsep
\arrayrulewidth
\begin{array}{rr rr rr rr r}
 1& 1& 1& 1& 1& 1& 1& 1& 1\\ 
 1&  &-1& 1&  &-1& 1&  &-1\\ 
 1&-1& 1& 1&-1& 1& 1&-1& 1\\ 
 1& 1& 1&  &  &  &-1&-1&-1\\ 
 1&  &-1&  &  &  &-1&  & 1\\ 
 1&-1& 1&  &  &  &-1& 1&-1\\ 
 1& 1& 1&-1&-1&-1& 1& 1& 1\\ 
 1&  &-1&-1&  & 1& 1&  &-1\\ 
 1&-1& 1&-1& 1&-1& 1&-1& 1   
\end{array}
\]

\end{example}


\section{Multiplicative inverse means mirroring}

For studying the groups of symmetries of the black and white image $I$ is enough to understand the 
symmetries for the fundamental block $F(p,m)$. All groups of symmetries we are looking for are subgroups
of the dihedral group of symmetries $D_8$ of the square. We start with the most non-symmetric case, the case 
of Pascal's Triangle:

\begin{lemma} \label{m=0}
If $m=0$ the group of symmetries consists of two elements: 
the identity and the reflection  through the first 
diagonal.
\end{lemma}

{\bf Proof}: In $F(p,0)$ for $0 \leq i,j \leq p-1$ :
\[a_{i,j} = 0 \,\,\leftrightarrow \,\, p \, | \, f(i,j) = \binom{i+j}{i} \,\,\leftrightarrow \,\, i+j \geq p.\]
So exactly the elements strictly below the second diagonal are $0$.\qed

Now a very small theory shall explain the other cases. 

\begin{definition}\rm
For a matrix $A$ we define the mirror-image $\Sigma A$ using the definition of a matrix 
as a family of column-vectors. If $A=(\vec a_1, \dots, \vec a_n)$ 
then $\Sigma A=(\vec a_n, \dots , \vec a_1)$.
\end{definition}

\begin{definition}\rm
For $m \neq 0$   we define the operator ${\cal O}$ acting over the fundamental block $F(p,m)$ in the following way:
\[ {\rm For}\,\,i=0 \,\,{\rm to}\,\, p-1, \,\,{\rm one}\,\,{\rm divides}\,\, {\rm  the} \,\,
{\rm row}\,\, i \,\, {\rm by}\,\, (-m)^i.\]
The result is denoted by ${\cal O} F(p,m)$.
\end{definition}

\begin{lemma}\label{trick}
For all finite fields $\F_q$ and for all $m \in \F_q \setminus \{0\}$ the following identity holds:
\[ {\cal O} F(p,m) = \Sigma F(p,m^{-1}).\]

\end{lemma}

{\bf Proof}: The Lemma follows from the following claims:

\begin{enumerate}
\item The first row and the last column of  ${\cal O} F(p,m)$ consist only in ones.
\item For every  connected $2 \times 2$ sub-block of ${\cal O} F(p,m)$:
\[
\begin{array}{cc}
A & B \\
C & D
\end{array}
\]
 is true that $C = m^{-1} B + A + D$.
\end{enumerate}

The first claim follows from Lemma \ref{lastrow} and from the definition of the operator ${\cal O}$: one divides
exactly with the elements of the last column.

We prove the second claim. Let $(a,\, b\,\,|\,\,c,\,d)$ be the corresponding elements in $F(p,m)$. They fulfill
the equality:
\[d=ma+b+c.\]
Using the definition of ${\cal O}F(p,m)$, we see that:
\[ A = \mu a, \,\,\,\, B = \mu b, \]\[ C = (-m)^{-1}\mu c, \,\,\,\, D = (-m)^{-1} \mu d,\]
where $\mu = (-m)^i$ for some $i$. It follows that:
\[ C = (-m)^{-1} \mu c = (-m)^{-1} \mu (d-ma-b) = (-m)^{-1} \mu d +  \mu a - (-m)^{-1} \mu b = D + A + m^{-1}B.\]
\qed

For the next Corollary recall from the proof of  \ref{fractalInGeneral} that $\delta F(p,m)$ is the matrix obtained
by substituting every element of $F(p,m)$ with $1$ if it is $\neq 0$. Also, recall that the elements of $F(p,m)$ are 
by $a_{i,j}$ and denote the elements of $F(p, m^{-1})$ with $a_{i,j}'$.

\begin{lemma}\label{k4inside} The following statements follow directly from Lemma \ref{trick}:
\begin{enumerate}

\item For all $0 \leq i,j \leq p-1$:
\[  a_{i, p-1-j} ' = a_{i,j} (-m)^{-i}.\]

\item For all $0 \leq i,j \leq p-1$:
\[ a_{i,j} (-m)^{-i} = a_{p-1-j, p-1-i} (-m)^{j + 1 - p}. \]

\item If $m \in \F_q\setminus \{0\}$ then:
\[ \delta F(p,m) = \delta \Sigma F(p, m^{-1}) = \Sigma \delta F(p,m^{-1}).\]

\item If $m \in \F_q \setminus \{0\}$  the matrix $\delta F(p,m)$ allows two diagonal symmetries; and so
 all its tensor powers.
 
\item Given $m \in \F_q \setminus \{0\}$ fixed, some matrix $M_d$ contains zeros if and only if $M_1 = F(p,m)$ contains
zeros. If this takes places, then 
\[\deg ( m / \F_p) \leq \frac{p-1}{2}.\]

\end{enumerate}
\end{lemma}

{\bf Proof}: 
\begin{enumerate}
\item This is nothing as Lemma \ref{trick} written element-wise.

\item This is the symmetry of $F(p,m^{-1})$ through its first diagonal: just write the elements of $F(p,m^{-1})$ 
using the first statement of the present
Lemma as functions of the row-number and the corresponding element of $F(p,m)$. Concretely one has:
\[  a_{i, p-1-j} ' = a_{i,j} (-m)^{-i},\]
as in the precedent statement, 
\[ a_{i,p-1-j} ' = a_{p-1-j,i}',\]
which expresses the symmetry of $F(p,m^{-1})$ through its first diagonal, and
\[ a_{p-1-j,i}' = a_{p-1-j, p-1 -i} (-m)^{j+1-p},\]
which is an other instance of the first statement. Apply the transitivity.

\item It follows from the first statement:
\[ a_{i,j} = 0 \,\, \leftrightarrow \,\, a_{i, p-1-j} ' = 0.\]

\item For the reflexion through the first diagonal it is nothing to prove, because the recurrent law 
is symmetric. The symmetry through the second diagonal follows from the second
statement of the present Lemma:
\[ a_{i,j} = 0 \,\, \leftrightarrow \,\, a_{p-1-j, p-1-i} = 0.\]

\item The existence of zeros in $M_d$ is equivalent with the existence of zeros in $M_1 = F(p,m)$ because
$M_d$ is a tensor product of Frobenius conjugates of $F(p,m)$. Recall that the set of zeros of the fundamental
block $F(p,m)$ is symmetric regarding both diagonals, so there will be a zero with both coordinates $i,j \leq (p-1)/2$.
But the value of $a_{i,j}$ is a polynomial over $\F_p$ of degree $\min(i,j)$ in $m$.

\end{enumerate}
\qed

\begin{example}\rm
The last condition occurring here implies the existence of relatively less values of $m$ generating non-trivial
self-similar sets in arbitrary finite fields. Look at the case $\F_{19^2} = \F_{361}$ seen as $\F_{19}[x]$ where
$x^2 + 1 = 0$. Encode the element $ax+b$ in the natural number $19a+b$. I do not mention both $m$ and $m^{-1}$ because
they produce mirrored carpets. Also, if $m$ has been already mentioned, I don't mention its Frobenius $m^{19}$, because
it produces the same carpet. So, up to Frobenius and multiplicative inverse, one has non-trivial self-similar
carpets over $\F_{361}$ if and only if $m$ is equal with one of the following $29$ elements: 
$0$, $1$, $2$, $3$, $4$, $6$, $7$, $8$, $9$, $14$, $19$, 
$21$, $35$, $47$, $52$, $53$, $56$, $63$, $69$, $76$, $78$, $88$, $92$, $102$, $130$, $136$, $137$, $148$, $168$.
Values of $m \in \F_{361}$ which are not itself, inverses of, or Frobenius 
of elements in this list generate however interesting
coloured images: without zeros, but still tensor products of Frobenius conjugates of their fundamental blocks.
\end{example}


\section{$\F_p$ as a field of self-similar carpets}


\subsection{Complete classification}

\begin{theo}\label{symm}
Let $p$ be a prime and $m \in \F_p$. Exactly one of the following situations arrises:
\begin{enumerate}
\item $m=0$. In this case $I$ is a self-similar Pascal Triangle, $I$ is only symmetric through the first diagonal,
and the group of symmetries of $I$ is isomorphic with $S_2$.
\item $m = \pm 1$. In this case $I$ is a full square (for $m=-1$) or a nontrivial self-similar set
(for $m=1$) and the group of symmetries of $I$ is the full dyhedral group $D_8$ of the square.
\item $p \geq 5$, $m \in \F_p \setminus \{-1,0,1\}$. In this case $I$ is a non-trivial self-similar set
 and the group of symmetries of $I$ is generated by
the reflexions through the diagonals of the square. This group is isomorphic with Klein's group $K_4$. 
\end{enumerate}
\end{theo}

{\bf Proof}: The case $m=0$ follows completely from Lemma \ref{m=0}. Let now $m \in \F_p \setminus \{0\}$,  let
$K$ be the group generated by the symmetries through the both diagonals (isomorphic with Klein's group $K_4$) and let
$G$ be the group of symmetries of $I$. From Lemma \ref{k4inside} follows that $K \leq G \leq D_8$. 
If $m=-1$ then $I$ is the full square and trivially $G = D_8$. If $m = 1$ than it follows from Lemma \ref{k4inside} that:
\[ \delta F(p,1)  = \Sigma \delta F(p,1),\]
because $1^{-1}=1$ so $G$ is strictly bigger than $K$ which already has $4$ elements, hence $G = D_8$. 

Conversely, suppose that $G = D_8$. We exclude the trivial case $m = -1$. 
From Lemma \ref{m=0} it follows that $m \neq 0$. From Lemma \ref{existencezero}
it follows that $F(p,m)$ has only a zero in the second row ($i=1$), which is 
$a_{1,k} = 0$ for $0 < k < p-1$ such that $k = -(m+1)^{-1}$ in $\F_p$. If $a_{1,k}$ is not the
central element of the row $i = 1$ then  there are two zeros in this row: $a_{1,k}$ and its miror image through $\Sigma$.
This would be in contradiction with Lemma \ref{existencezero}. 
It follows that $ -(m+1)^{-1} = (p-1) /2 $ in $\F_p$, so $(m+1)^{-1} = 2^{-1}$, so $m=1$. 
\qed

\begin{cor}\label{cor}
If $p > 3$ and $m \in \F_p \setminus \{ -1 \}$ there are at least two zeros in $F(p,m)$.
\end{cor}

{\bf Proof}: In fact one can proof a litle bit more. If $m \in \{-2, -2^{-1}\}$ then the unique zero of the second row
$i=1$ lies on the intersection of this row with one of the diagonals, so its orbit under the action of $G$ has two elements.
If $m\in \F_p \setminus \{-2, -2^{-1}, -1, 0 \}$ then the orbit has four elements. For the case $m=0$ the existence of
many zeros is trivial, for $m=1$ a remark in the next section assure the existence of much more zeros than four.\qed


\subsection{The special case $m\in \{ -2, - 2^{-1} \}$: Diagonal Carpets}

For $m=-2$ one has $a_{1,1} = 0$. Mirror-symmetric: for $m=-2^{-1}$ one has $a_{1,p-2} = 0$. 
In fact, in these cases,  all the elements of odd index on the corresponding diagonal are zero!

\begin{definition}\rm Call {\bf first odd diagonal} (respectively {\bf second odd diagonal}) the following set of
indices:
\[D^+ = \{ (i,i) \,\,|\,\, 0 < i < p-1\,\, \wedge \,\, 2 \not | \, i\}.\]
\[D^- = \{ (i,j) \,\,|\,\, i + j = p-1\,\, \wedge \,\, 2 \not | \, i\}.\]

\end{definition}

\begin{theo}\label{odddiagonals} Let $p \geq 5$ be a prime. Following statements hold:
\begin{enumerate}

\item $D^+$ consists of zeros of $F(p, -2)$.

\item $D^-$ consists of zeros of $F(p, -2^{-1})$.

\item Moreover, the even elements on the respective diagonals are  $\neq 0$.
\end{enumerate}
\end{theo}

{\bf Proof}:

\begin{enumerate}

\item  We prove that for $m=-2\in \Z$ the recurrent function $f:\N \times \N 
\rightarrow \Z$ defined in the second section has the property $f(2s+1, 2s+1) = 0$ for all $k \in \N$. This follows
from the following identity:
\[
\sum \limits _{a=0} ^n (-2)^a \binom{n}{a} \binom{2n-a}{n-a} =
\begin{cases}
(-1)^s \binom{2s}{s}, & \,\,{\rm if}\,\, n=2s, \\
\,\,\,\,\,\,\,\,\,\,0, & \,\,{\rm if}\,\, n = 2s+1. 
\end{cases}
\]
This identity can be proved with Zeilberger's Algorithm, see \cite{PWZ} and \cite{K}. 
In fact, after running the software from \cite{K}, 
one gets the recurrent formula:
\[4(n+1)S(n) + (n+2)S(n+2) = 0,\]
where $S(n)$ is the sum on the left side of the equality. Starting with $S(0)=1$ and $S(1)=0$ one gets the result by
induction.  The author thanks Prof. Dr. Wolfram Koepf for kindly running his Maple package 
"Hypergeometric Summation" at author's request. Please note that this identity
is not the Reed - Dawson Identity, although a similar one.

\item This follows from the case $m = -2$ and the dualism from Lemma \ref{trick}. 
Note that the corresponding values of $f(n,k)$ are no more $0$ in $\Z$
but become $0$ in $\F_p$. 

\item This follows from the fact that the even diagonal values of $f(n.n)$ are $\neq 0$ in $\F_p$.
\end{enumerate}
\qed

We note that those regular diagonal zero are not the only zeros in general: starting with $p = 11$ there are a lot of 
sporadic (non-regular) zeros for $m \in \{-2, -2^{-1}\}$. One can now also prove a slight improvement of \ref{cor}:

\begin{cor}\label{improvcor} Let $m \neq -1$.
For $p=7$ there are at least three zeros in $F(p,m)$, and for $p \geq 11$ there are at least four zeros in $F(p,m)$.
\end{cor}

{\bf Proof}: The only remained problem was $m \in \{-2, -2^{-1}\}$, which is now trivial 
applying Theorem \ref{odddiagonals}. In fact, there are many more sporadic (non-diagonal) 
zeros for these values of $m$. \qed


\subsection{The special case $m=1$: Cross-carpets}

The only one fully symmetric non-trivial case (where $p$ is odd and $m=1$) is worth for a closer look. This is exactly
the case of the spectacular Passoja-Lakhtakia Carpets, described in \cite{PL}. The author noticed after working some
weeks with these carpets, that the infinite symmetric matrix $(f(n,k))$ for $m=1$ is known as the double sequence of the
Delannoy Numbers, and is the answer to the following combinatoric problem: {\it In how many ways one can go from the
left-highest corner of a matrix to the element $a_{n,k}$ making only South, South-East and East steps.} For a brief
history of the Delanoy numbers, see \cite{BS}, also in the On-line Encyclopedia of Integer Sequences.

\begin{definition}\rm Let us call $N = \{(i,j) \,|\,a_{i,j} =0\}$ the set of {\bf zeros} of the fundamental cell $F(p,1)$.
The set:
\[C = \{(\frac{p-1}{2},i)\,; \,(i,\frac{p-1}{2})|\, 0 \leq i \leq
p-1 \,\, \wedge \,\, 2\not | \,i\}\]
shall be called the {\bf cross}, 
and $S = N \setminus C$  the set of {\bf sporadic} zeros. We call the elements of the cross {\bf regular zeros}. The
next Corollary proves that the elements of the Cross are really zeros of $F(p,1)$.
\end{definition}

\begin{cor}\label{PL} If $p$ is an odd prime, the fundamental block $F(p,1)$ has the following properties:
\begin{enumerate}

\item For all $ 0 \leq k  \leq p-1$:
\[a_{p-1, k} = (-1)^k.\]

\item For all  $n$ and $k$  with  $0 \leq n,k \leq p-1$,
\[a_{n,k} = (-1)^n a_{n,p-1-k}.\]

\item The Cross $C$ consists of zeros of $F(p,1)$.

\end{enumerate}
\end{cor}

{\bf Proof}:
\begin{enumerate}
\item This is exactly Lemma \ref{lastrow}.

\item According to Lemma \ref{trick},
\[ {\cal O} F(p,1) = \Sigma F(p,1).\]
If $k = 2s$, ${\cal O}$ operates by multiplication with $1$, so the even rows are centrally symmetric. If $k=2s+1$,
${\cal O}$ operates by multiplication with $-1$, so odd rows are antisymmetric.

\item This follows easily from the last statement because for $k$ odd, $p$ being also odd: 
\[a_{k,\frac{p-1}{2}} = {(-1)}^k a_{k,
(p-1)-\frac{p-1}{2}} = - a_{k,\frac{p-1}{2}},\] implies $a_{k,\frac{p-2}{2}} = 0$. 
Apply now the symmetry of $\delta F(p,1)$.

\end{enumerate}
\qed

\begin{example}\rm
Here one sees only the bord and the zeros of the fundamental block $F(13,1)$:

\[
\arraycolsep
\arrayrulewidth
\begin{array}{ccccccccccccc}
+1&+1&+1&+1&+1&+1&+1&+1&+1&+1&+1&+1&+1\\ 
+1&\cdot &\cdot &\cdot &\cdot &\cdot &0&\cdot &\cdot &\cdot &\cdot &\cdot &-1\\ 
+1&\cdot &0&\cdot &\cdot &\cdot &\cdot &\cdot &\cdot &\cdot &0&\cdot & +1\\ 
+1&\cdot &\cdot &\cdot &\cdot &\cdot &0&\cdot &\cdot &\cdot &\cdot &\cdot &-1\\ 
+1&\cdot &\cdot &\cdot &\cdot &\cdot &\cdot &\cdot &\cdot &\cdot &\cdot &\cdot & +1\\ 
+1&\cdot &\cdot &\cdot &\cdot &\cdot &0&\cdot &\cdot &\cdot &\cdot &\cdot &-1\\ 
+1&0&\cdot &0&\cdot &0&\cdot &0&\cdot &0&\cdot &0& +1\\ 
+1&\cdot &\cdot &\cdot &\cdot &\cdot &0&\cdot &\cdot &\cdot &\cdot &\cdot &-1\\ 
+1&\cdot &\cdot &\cdot &\cdot &\cdot &\cdot &\cdot &\cdot &\cdot &\cdot &\cdot & +1\\ 
+1&\cdot &\cdot &\cdot &\cdot &\cdot &0&\cdot &\cdot &\cdot &\cdot &\cdot &-1\\ 
+1&\cdot &0&\cdot &\cdot &\cdot &\cdot &\cdot &\cdot &\cdot &0&\cdot & +1\\ 
+1&\cdot &\cdot &\cdot &\cdot &\cdot &0&\cdot &\cdot &\cdot &\cdot &\cdot &-1\\ 
+1&-1&+1&-1&+1&-1&+1&-1&+1&-1&+1&-1&+1 
\end{array}
\]
\end{example}

The primes $3$, $5$, $7$, $11$, $19$ have only regular zeros in $F(p,1)$. 
$13$ is the first odd prime with sporadic zeros, followed by $17$. 
By all other primes tryed out by the author (from $23$ to $599$) 
there are lots of sporadic zeros in the fundamental block $F(p,1)$.


\section{Aperiodic tilings}

An interesting implication of the Theorem \ref{tensorpower} is that we get for free the existence of the aperiodic
tilings for the quarter of plane. In order to preserve the matrix-friendly notation, we will refer here to the quarter
$\R_+ \times \R_-$. We use a definition suggested in \cite{Sen}:

\begin{definition}\rm
Let $A \subseteq \R^2$ be an unbounded subset.

A finite set of bounded polygons is called {\bf aperiodic tiling} of $A$ if both conditions are fulfilled:

(1) $A$ can be partitioned in copies of the polygons.

(2) There is no translation of the partition carrying every polygon in a copy of itself.

\end{definition}

\begin{cor}
Let $\F_q$ be some finite field of characteristic $p$ and $m \in \F_q \setminus \{ -1\}$ 
such that the fundamental block $F(p,m)$ contains at 
least one zero.  The self-similar carpet defined by $m$ over $\F_q$ can be seen as aperiodic tiling
of $\R_+ \times \R_-$, done with a finite set of polygons.
\end{cor}

{\bf Proof}: 

{\bf Case} 1: $m=0$. If $m = 0$ one has a set of  $\leq p^2$ coloured square tiles. 
The set of colours is in a bijection with the elements of the prime  field  $\F_p$, as before the colour $0$
is called white. The tiles are coloured as follows:


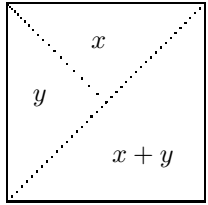
\begin{picture}(120,120)(-180,0)

\centering


\qbezier(0,0)(25,0)(75,0)
\qbezier(75,0)(75,50)(75,75)
\qbezier(75,75)(50,75)(0,75)
\qbezier(0,75)(0,50)(0,0)


\qbezier[40](0,0)(38,38)(75,75)
\qbezier[20](0,75)(5,70)(37.5,37.5)


\put(32,58){$x$}
\put(10,38){$y$}
\put(40,15){$x+y$}

\end{picture}

\vspace*{1cm}


We call the three teritories of a tile North, West and Big South-East. This set of tiles works with the following rules:
\begin{enumerate}
\item They can be only translated, but not rotated. 
\item The edges touching the positive part of the axis $Ox$ are nordic edges in tiles with Nord  coloured with $1$.
\item The edges touching the negative part of the axis $Oy$ are western edges in tiles with West coloured with $1$.
\item If two tiles are touching along an edge, this can be only a Nord - South or a West - East contact along a whole edge.
\item In the case of a Nord - South contact, the North of the south tile has the same colour as the Big South-East of 
the north tile.
\item In the case of a East - West contact, the West of the eastern tile has the same colour than the Big South-East 
of the western tile.
\end{enumerate}
The rules implies that the tile touching $(0,0)$ has the colours West $=1$, North $=1$ and Big South-East $=2$ in $\F_p$.


{\bf Case} 2: $m \neq 0$. If $m \neq 0$ one has a set of $\leq r^3 +2$ coloured tiles, 
where $r$ is the number of elements of the 
field $\F_p[m]$. Like in the case $m = 0$ to different elements correspond different colours. To $0$ corresponds $white$.

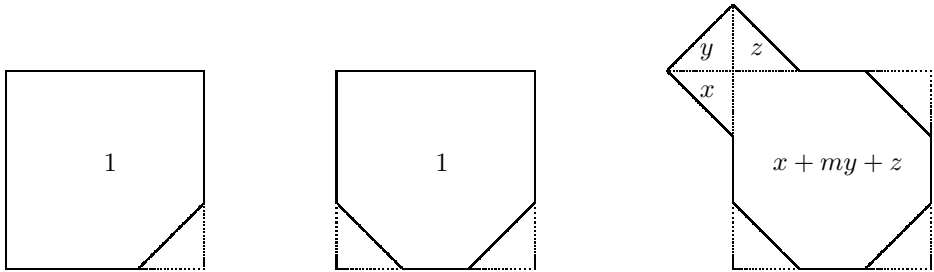
\begin{picture}(400,120)(-50,0)

\centering


\qbezier(0,0)(25,0)(50,0)
\qbezier(50,0)(62.5,12.5)(75,25)
\qbezier(75,25)(75,50)(75,75) 
\qbezier(0,75)(25,75)(75,75)
\qbezier(0,0)(0,25)(0,75)


\qbezier[20](50,0)(55,0)(75,0)
\qbezier[20](75,0)(75,10)(75,25)


\qbezier(150,0)(160,0)(175,0)
\qbezier(175,0)(187.5, 12.5)(200,25)
\qbezier(200,25)(200,50)(200,75)
\qbezier(125,75)(150,75)(200,75)
\qbezier(125,25)(125,50)(125,75)
\qbezier(125,25)(137.5, 12.5)(150,0)


\qbezier[20](125,0)(130,0)(150,0)
\qbezier[20](125,0)(125,5)(125,25)

\qbezier[20](175,0)(180,0)(200,0)
\qbezier[20](200,0)(200,5)(200,25)


\qbezier(300,0)(312.5, 0)(325,0)
\qbezier(325,0)(337.5, 12.5)(350,25)
\qbezier(350,25)(350,30)(350,50)
\qbezier(350,50)(337.5, 62.5)(325,75)
\qbezier(325,75)(320,75)(300,75)
\qbezier(300,75)(287.5,87.5)(275,100)
\qbezier(275,100)(262.5,87.5)(250,75)
\qbezier(250,75)(255,70)(275,50)
\qbezier(275,50)(275,45)(275,25)
\qbezier(275,25)(280,20)(300,0)


\qbezier[40](275,50)(275,75)(275,100)
\qbezier[40](250,75)(275,75)(300,75)

\qbezier[20](275,0)(280,0)(300,0)
\qbezier[20](275,0)(275,5)(275,25)

\qbezier[20](325,0)(330,0)(350,0)
\qbezier[20](350,0)(350,5)(350,25)

\qbezier[20](350,50)(350,55)(350,75)
\qbezier[20](325,75)(330,75)(350,75)


\put(37,37){$1$}
\put(162.5,37){$1$}
\put(262.5, 65){$x$}
\put(262.5, 81){$y$}
\put(281.5, 81){$z$}
\put(290, 37.5){$x+my+z$}

\end{picture}

\vspace*{1cm}

These are the tiles of types one, two and three, in this order from left to the right. 
The four regions of the tiles of type three are called West (marked with $x$), North-West (marked with $y$), North 
(marked with $z$) and Big Body. North, North-West and West form the so-called Active Side. 
This set of tiles works with the following rules:
\begin{enumerate}
\item Only one tile of type one shall be used. This tile is the only one touching $(0,0)$.
\item Along the axis $Ox$ one puts tiles of type two, with the north edge on the axis. Along the axis $Oy$ one puts 
tiles of type two rotated with $90^\circ$ to their left, showing inside the quarter of plane.
\item All other tiles used are tiles of type three. They might be translated, but never rotated.
\item The three regions of an Active Side (North, North-West and West) always have the same colours as the Big Body 
that they are respectively touching: North has the same 
colour as the Big Body of the northern neighbor, North-West has the 
same colour as the Big Body of the north-west neighbor and West has 
the same colour as the Big Body of the 
western neighbor. Other said, one always combines tiles to get a composite covering of monochromatic squares.
\end{enumerate}

It follows from the Theorem \ref{tensorpower} that the first set of tiles constructs the convergents to Pascal's
Triangle in $\F_p$ and that the second set of tiles constructs the convergents to the self-similar carpet defined
by $m \in \F_q$. Why are these tilings aperiodic? Suppose that one of them is invariant to some translation with a 
vector $ \vec v \in \R^2$ of length $v$. In both cases there are arbitrary big bounded white areas, 
so there is a white area containing a disk of radius $> 2v$. This disk shall intersect its own translation, so the 
bouded white area also. This is a contradiction. \qed

\begin{remark}\rm
Most of the conditions can be naturally encoded in the tiles by modifying their form using local convex versus 
concave adds.
\end{remark}


\section{Commentaries}

\begin{enumerate}

\item Can we better understand the sporadic zeros for $\F_p$ and $m=1$? The Hausdorff (fractal) dimension 
depends only of the total number of zeros in the fundamental block. The same question for $\F_q$ and arbitrary $m \in \F_q$.

\item Is it true that two zeros of the fundamental block cannot have a common edge? The fundamental 
block contains sometimes edge-neighbors
with equal value: try for example $\F_{11}$ and $m=1$ where $a_{4,2}=a_{4,3}=a_{3,3}=a_{3,4}= 8$. 
The author found some cases of
zeros with common vertex in a fundamental block, but no example with common edge.

\item For the case $\F_p$, $p$ odd prime and $m=1$ the symmetries proved in Corollary \ref{PL} 
justifies the following strategy of
representation: choose {\bf white} for $k=0$ and a list of colours such that
\[ \forall \, k \in \F_p \,\, {\rm color} (k) = {\rm color}(p-k). \]

\item For some integer $n$ with prime-decomposition  $n=p_1^{k_1} \dots p_s^{k_s}$, 
the ring $\Z / n \Z$ is isomorphic with the product of finite rings
$\Z / p_1^{k_1} \Z \times \dots \times \Z / p_s^{k_s} \Z$. From this reason, 
carpets over $\Z_n$  are overlappings of carpets modulo $p^k$. 
Can we understand the (black and white, or coloured) carpets modulo prime-powers? 
In the experiments made by the author  very complicated patterns arose. For example, the black and white matrix
$\delta F(p^2, 1)$ presents a cross that is periodically interupted by black and white patterns $\delta F(p,m)$.
Such patterns appear also in positions corresponding to the sporadic zeros of $F(p,1)$, with coordinates multiplied with
$p$. A lot of new sporadic zeros can be remarked in $F(p^2, 1)$. 

\item Over fields $\F_p$ the rule $a_{n,k} = x a_{n-1,k-1}+ ya_{n,k-1}+za_{n-1,k}$ generates very complicated structures.
At depth $\geq 3$ one recognizes the zeros from the case discussed here, with the following sensible difference: in the 
holes of the carpets there are different periodic structures, instead of the constant zero hole. 
I call this structures {\bf patchwork carpets}.
For exceptional values
of $x$, $y$ und $z$ one gets one of the already known self-similar carpets, or periodic structures, or constant structures.

\end{enumerate}


\end{document}